\documentclass[a4paper,12pt]{article}
 \usepackage[top=2.5cm,bottom=2.5cm,left=2.5cm,right=2.5cm]{geometry}
 \usepackage{cite, amsmath, amssymb}
\newtheorem{prop}{Proposition}

 \title{\vspace*{2.5cm}{\bf Kemeny's constant and  the Kirchhoff index for the cluster of highly symmetric graphs} }
\date{}
\author{{\bf Jos\'e Luis Palacios}\\
{\it Electrical and Computer Engineering Department,}\\{\it The University of New Mexico, Albuquerque, NM 87131, USA}\\
jpalacios@unm.edu\\
and\\
{\bf Greg Markowsky}\\
{\it Department of Mathematics,}\\
{\it Monash University, Melbourne, Australia}\\greg.markowsky@monash.edu}

\begin{document}
 \thispagestyle{empty} 
\maketitle

%\begin{center}(Received December 11, 2015)
%\end{center}

\begin{abstract} 
We find closed form formulas for the Kemeny's constant and the Kirchhoff index for the cluster $G_1\{G_2\}$ of two highly symmetric graphs $G_1$, $G_2$, in terms of the parameters of the original graphs. We also discuss some necessary conditions for a graph to be highly symmetric. \end{abstract}

\section{Introduction}

Let $G=(V,E)$ be a finite simple connected graph with vertex set $V=\{1, 2, \ldots, n\}$, edge set $E$ and  degrees $\Delta =d_1\ge d_2\ge \cdots \ge d_n=\delta$. 
The simple random walk on $G$ is the Markov chain $X_n, n\ge 1$ that jumps from one vertex of $G$ to a neighboring vertex with uniform probabilities.  If $P$ is the transition matrix of this chain, the stationary distribution of the random walk is the unique probabilistic vector $\pi$ that satisfies  $\pi P=\pi$ and that can be explicitly given as $\pi_=\frac{d_i}{2|E|}$. The hitting time $T_b$ of the vertex $b$ is defined as $T_b=\inf \{n: X_n=b\}$ and its expectation, when the process starts at vertex $a$ is denoted by $E_aT_b$.  Sometimes we will use superscripts $G, G_1$, as in $E_aT^G_b$ and $E_aT^{G_1}_b$ to distinguish between the walk on the graph $G$ or a subgraph $G_1$.

Given two graphs $G_1=(V_1, E_1)$, $G_2=(V_2, E_2)$, the cluster $G_1\{G_2\}$ (following the notation of \cite{ZYL}) is the graph obtained by selecting one vertex of $G_2$ and identifying this vertex with every vertex of $G_1$, so that $|V_1|$ copies of $G_2$ are glued to $G_1$, one for each vertex. For example, if $P_2$ is the 2-path, $G_1\{P_2\}$ consists of adding to every vertex of $G_1$ another vertex linked with a single edge.

The Kemeny constant is defined as 
$$K=\sum_j \pi_j E_iT_j,$$ 
which is a constant independent of $i$. We direct the reader to \cite{GS} for this fact and all other probabilistic details. For graph notions we refer to \cite{BH}.

The Kirchhoff index is a molecular descriptor defined on an undirected connected graph as
$$R(G)=\sum_{i<j} R_{ij},$$
where $R_{ij}$ is the effective resistance between $i$ and $j$ computed on the graph when it is thought of as an electric network with unit resistors on each edge. Kirchhoff indices and electric resistances on graphs are important notions in mathematical chemistry, and have attracted attention from the pure mathematics community as well. Recent works related to the results of this paper include \cite{atik, ciardo, deville, Faught, geng, li, lei, zhou, zhou2}.

Chandra et al showed in \cite{C} that there is a close relationship between hitting times and effective resistances
\begin{equation}
\label{chandra}
E_aT_b+E_bT_a =2|E|R_{ab}.
\end{equation}

In \cite{P1}, \cite{PR}, it was noticed that for $d$-regular graphs there is a simple relationship between the Kirchhoff index and the Kemeny constant:
\begin{equation}
\label{regular}
R(G)=\frac{|V|}{d}K,
\end{equation}
but the relationship between $K$ and $R(G)$ is not as straightforward when the graph is not regular. See \cite{W}, \cite{X}, \cite{KD}, \cite{PM} for calculations of $K$ involving non-regular graphs and their Kirchhoffian indices.

In a series of articles (\cite{P1}, \cite{PR2010}, \cite{PR}, \cite{BCPT}) we devoted to the Kirchhoff index, the Kemeny constant took a back seat, but we can highlight a number of results for $K$ collected in those articles. In \cite{PR2010}, writing $K$ in terms of the eigenvalues of the transition matrix of the random walk on $G$ as
\begin{equation}
\label{keigen}
K=\sum_{j=2}^n\frac{1}{1-\lambda_j},
\end{equation}
we found that 
$$K\ge \frac{(n-1)^2}{n},$$
and the lower bound is attained by the complete graph $K_n$.

In case the graph is bipartite we can improve slightly the bound to
$$K\ge \frac{2n-3}{2},$$
and the lower bound is attained by the complete bipartite graph $K_{\frac {n}{2}, \frac{n}{2}}$. The maximum value of $K$ over all graphs is not known, to the best of our knowledge, though we can say that it is at least cubic in $n$, since for the $(1/3, 1/3, 1/3)$-barbell graph which consists of two copies of the complete $K_{n/3}$ graph attached at the endpoints of a path of length $n/3$, it is known (see \cite{PR}) that $\lambda_2\ge 1-\frac{c}{n^3}$, where $c=54+O(1/n)$, and therefore by (\ref{keigen}) we have that for this barbell graph
$$K\ge \frac{n^3}{c}.$$
In \cite{BCPT}, using the inequalities
\begin{equation}
\label{notsosimple}
\frac{n}{\Delta}K \le R(G)\le \frac{n}{\delta} K
\end{equation}
that generalize (\ref{simple}), we found several bounds for the Kirchhoff index, through bounds for $K$ that were obtained with majorization, namely,  for any graph $G$ we have 
$$
K\ge \frac{1}{1+\frac{\sigma}{\sqrt{n-1}}}+\frac{(n-2)^2}{n-1-\frac{\sigma}{\sqrt{n-1}}},
$$
where 
$$\sigma^2=\frac{1+\sum_{i=2}^n\lambda_i^2}{n},$$
and
$$K\le \frac{n-k-2}{1-\lambda_2}+\frac{k}{2}+\frac{1}{\theta},$$
where $\displaystyle k=\lfloor \frac{\lambda_2(n-1)+1}{\lambda_2+1}\rfloor$ and $\theta=\lambda_2(n-k-2)-k+2$. 

Several other more specific bounds for $K$ were found in \cite{BCPT} that the reader could check, for instance,  if $G$ is $d$-regular with diameter $D$ then
$$K\ge \frac{1}{1+\frac{2D}{d(D+1)}}+\frac{(n-2)^2}{n-1-\frac{2D}{d(D+1)}}.$$

Some recent articles have dealt with the Kemeny constant of some composite graphs, such as the subdivision and triangulation of a graph (see \cite{X1} \cite{X2}) or barbell graphs (see \cite{Br}).

In the previous recent article \cite{PM} we found explicit relationships between $K$, $R(G)$ and $R^*(G)$ for a family of non-regular graphs (some of which are barbell graphs) created from conjoining several copies of a ``highly symmetric" (HS) graph for which
\begin{equation}
\label{vaya}
E_aT_b = E_bT_a,
\end{equation}
for all $a, b \in G$.  The family of walk-regular graphs satisfies (5).  A graph is walk-regular if the number of $k$-long walks, $k\ge 2$,  starting and ending at a vertex $v$ is the same for all $v\in V$.  This family contains the families of vertex-transitive, regular edge-transitive and distance regular graphs.  A nice feature of these HS graphs,  that simplifies computations with hitting times,  is the fact that (\ref{vaya}) and (\ref{chandra}) imply
\begin{equation}
\label{simple}
E_aT_b=|E|R_{ab}.
\end{equation}
Another important feature of these HS graphs is that the computation of the Kirchhoff index is also simplified, since if we define $R(i)=\sum_j R_{ij}$ then it is clear that for any graph
$R(G)=\frac{1}{2}\sum_i R(i)$, but for walk-regular graphs we have that $R(i)=\frac{2}{nd}\sum_j |E| R_{ij}=\frac{2}{d} \sum_j \pi_j E_iT_j=\frac{2}{d}K$, so that $R(i)$ is independent of $i$ and 
\begin{equation}
\label{easy}
R(G)=\frac{n}{2}R(i),
\end{equation}
for any $i\in V$.

In this article we continue exploring the relationships between $K$ and $R(G)$ for another family of non-regular graphs,  the cluster $G=G_1\{G_2\}$, where both $G_1$ and $G_2$ are HS, and express $K$ in terms of $K_1$ and $K_2$, $R(G)$ in terms of  $R(G_1)$ and $R(G_2)$  and  $R(G)$ in terms of $K$ When $G=G_1\{G_1\}$, yielding an expression more involved than (\ref{regular}).    The technique used in this article for finding $K$ in terms of $K_1$ and $K_2$ is similar to that in \cite{PM}, finding all hitting times of the random walk on $G$ in terms of effective resistances,  which is facilitated by the fact that (\ref{simple}) and (\ref{easy}) hold in $G_1$. This differs from the bulk of articles in the literature on Kemeny's constant, that rely on the analysis of eigenvalues.  We end the article with a discussion on some necessary conditions for a graph to be highly symmetric.

\section{The cluster $G_1\{G_2\}$}

We start with  two HS  graphs $G_1=(V_i, E_i)$, with $|V_i|=n_i$,  degrees $d_i$, and Kemeny constant $K_i$, $1\le i \le 2$, and consider its cluster $G_1\{G_2\}=G=(V, E)$.  Then it is clear that $|V|= n_1n_2$, $|E|=|E_1|+n_1|E_2|$ and the degrees of $G$ are either $d_1+d_2$ (for the vertices of the original $G_1$ or $d_2$ (for all others).   Under these conditions we have
\begin{prop}
For any HS graphs $G_1$, $G_2$ and their cluster $G$ defined above we have
\begin{equation}
\label{Bien}
K=\frac{|E|}{|E_1|}K_1+\frac{n_1(2|E|-|E_1|)}{|E|}K_2.
\end{equation}
\end{prop}

{\bf Proof.} We compute the summation of all expected values of hitting times, normalized with the stationary distribution,  when started at a vertex $c\in V_1$.  Let us denote by $d$ any point in a copy of $G_2$ glued to $G_1$ through $c$, we notice that
$$E_cT^G_ d+E_dT^G_c=2|E|R_{cd}$$
but $E_dT_c^G=E_dT_c^{G_2}=|E_2|R_{cd}$ and solving we get
\begin{equation}
\label{ayuda}
E_cT^G_d=(2|E|-|E_2|)R_{cd}.
\end{equation}
Now let $b\in V_1$, since (\ref{vaya}) holds in $G_1$, it is clear that
$$E_cT^G_b=|E|R_{cb}=(|E_1|+n_1|E_2|)R_{cb}=E_cT^{G_1}_b+n_1|E_2|R_{cb}.$$
Therefore
$$\sum_{b\in V_1} \pi_bE_cT^G_b=\frac{d_1+d_2}{2|E|}\sum_{b\in V_1}(E_cT^{G_1}_b+n_1|E_2|R_{cb})$$
$$=\frac{d_1+d_2}{2|E|}\frac{2|E_1|}{d_1}\sum_{b\in V_1}\frac{d_1}{2|E_1|}E_cT^{G_1}_b+\frac{d_1+d_2}{2|E|}n_1|E_2|\sum_{b\in V_1}R_{cb}$$
$$=\frac{|E_1|(d_1+d_2)}{|E|d_1}K_1+\frac{d_1+d_2}{2|E|}n_1|E_2|R^1(c),$$
where $R^1(c)$ denotes the sum of effective resistances from $c$ just over $G_1$, and this is
$$=\frac{|E_1|(d_1+d_2)}{|E|d_1}K_1+\frac{(d_1+d_2)|E_2|}{|E|}R(G_1),$$
where in the last equality we have used (\ref{easy}), and this in turn happens to be
\begin{equation}
\label{primera}
=\frac{|E_1|(d_1+d_2)}{|E|d_1}K_1+\frac{n_1|E_2|(d_1+d_2)}{|E|d_1}K_1=\frac{d_1+d_2}{d_1}K_1.
\end{equation}

Also, if $b\in V-V_1$ we have
$$E_cT^G_b=E_cT^G_{g(b)}+E_{g(b)}T^G_b,$$
where $g(b)\in V_1$ is the point of contact between $G_1$ and the copy of $G_2$ where $b$ lies, so using (\ref{ayuda}) this becomes
$$E_cT^G_b=E_cT^G_{g(b)}+(2|E|-|E_2|)R_{g(b)b},$$
and we can write:
$$\sum_{b\in V-V_1} \pi_bE_cT^G_b=\frac{d_2}{2|E|}\sum_{b\in V-V_1}E_cT^G_b$$
\begin{equation}
\label{intermediate}
=\frac{d_2}{2|E|}\sum_{b\in V-V_1}E_cT^G_{g(b)}+\frac{d_2(2|E|-|E_2|)}{2|E|}\sum_{b\in V-V_1}R_{g(b)b}.
\end{equation}
In the first sum in (\ref{intermediate}), to visit all the states in a copy of $G_2$ the random walk must visit $g(b)$, the glueing point,  $n_2-1$ times, so after relabeling the vertices, the sum becomes
\begin{equation}
\label{medio}
\frac{(n_2-1) d_2}{2|E|}\sum_{b\in V_1} E_cT_b^G=\frac{(n_2-1)d_2}{d_1+d_2}\sum_{b\in V_1}\pi_b E_cT_b=\frac{(n_2-1)d_2K_1}{d_1},
\end{equation}
where in the last equality we have used (\ref{primera}). The second sum in (\ref{intermediate}), which runs through the $n_1$ copies of $G_2$, with the help of (\ref{easy})  can be written as
$$\frac{n_1d_2(2|E|-|E_2|)}{2|E|}R(g(b))=\frac{n_1d_2(2|E|-|E_2|)}{n_2|E|}R(G_2),$$
and using (\ref{regular}) this becomes
\begin{equation}
\label{segunda}
\frac{n_1(2|E|-|E_1|)}{|E|}K_2.
\end{equation}

Putting (\ref{primera}) (\ref{medio}) and (\ref{segunda}) together we get
$$K=\sum_{b\in V}\pi_bE_cT_b=\sum_{b\in V_1}\pi_bE_cT_b+\sum_{b\in V-V_1}\pi_bE_cT_b$$
$$=\frac{d_1+n_2d_2}{d_1}K_1+\frac{n_1(2|E|-|E_1|)}{|E|}K_2$$
$$\frac{|E|}{|E_1|}K_1+\frac{n_1(2|E|-|E_1|)}{|E|}K_2~~\bullet$$
\vskip .2 in
Regarding the Kirchhoff index we can prove the following generalization of a result in \cite{ZYL}, where their hypothesis is that $G_2$ be vertex-transitive.  We generalize it to $G_2$ being HS. 
\begin{prop}
For arbitrary $G_1$ and $G_2$ HS and its cluster $G$  we have
\begin{equation}
\label{Bueno}
R(G)=n_2^2R(G_1)+(2n_1^2-n_1)R(G_2).
\end{equation}
\end{prop}

{\bf Proof}. We imitate the proof of this result in \cite{ZYL}, where they require that the quantities $R(i)$ be constant in $G_2$, and as discussed above this happens when $G_2$ is HS $\bullet$
\vskip .2 in
As an immediate corollary of the two previous propositions we get the following
\begin{prop}
For any $n$-vertex $d$-regular HS $G_1$ and its cluster $G=G_1\{G_1\}$  we have
\begin{equation}
\label{gorda1}
K=(3n+1)K_1,
\end{equation}
\begin{equation}
\label{gorda2}
R(G)=n(3n-1)R(G_1),
\end{equation}
and
\begin{equation}
\label{gorda3}
R(G)=\frac{n^2(3n-1)}{(3n+1)d}K.
\end{equation}
\end{prop}

\section{Necessary conditions for graphs to be highly symmetric}

The importance of symmetric hitting times for our results has led us to examine the question of when a graph can be highly symmetric. It is asserted in \cite{george1} that the family of walk-regular graphs is highly symmetric. A graph is walk-regular if the number of $k$-long walks, $k\ge 2$,  starting and ending at a vertex $v$, is the same for all $v\in V$. Every such graph must be regular, and this family contains the families of vertex-transitive and distance regular graphs. However, in practice, it may be difficult to determine whether a graph is walk-regular. Furthermore, the condition of walk-regularity is not known to be necessary for a graph to be highly symmetric. We have been interested in necessary conditions that are easily checked, and have found several such (though clearly not sufficient) conditions for a graph to be highly symmetric. First, a definition: we will say that a vertex $i$ in a graph is {\it resistance regular} if $R_{ij}=R_{ik}$ for any $j,k$ adjacent to $i$.

\begin{prop}

\begin{itemize}
    \item[$(i)$] If $G$ is a highly symmetric graph, and every vertex of $G$ is resistance regular, then $G$ is regular.
    
    \item[$(ii)$] If $G$ is a highly symmetric graph, then $G$ cannot have a resistance-regular cut-vertex.
    
    \item[$(iii)$] If $G$ is a highly symmetric graph with a cut-edge, then the two components that result from cutting the edge must have the same number of edges.

    \item[$(iv)$] A highly symmetric graph cannot contain two cut-edges.
\end{itemize}

\end{prop}

Before proving this, we makes several remarks. $(i)$ proves that no tree is highly symmetric other than the trivial one on two vertices, which was already known (\cite{george2}). It also shows that a highly symmetric edge-transitive graph must be regular; this was also shown in \cite{PM}, where it is further proved that a regular edge-transitive graph is automatically highly symmetric. A simple example of a graph to which $(ii)$ applies would be two $n$-gons with a single vertex from each conjoined, and other examples include the friendship graphs. The conditions in $(iii)$ and $(iv)$ are self-explanatory

{\bf Proof:} $(i)$ Since $G$ is connected, our conditions imply that the resistance across each edge is a constant. Therefore the hitting time across each edge is a constant $C$, and it follows from this that the return time $E_i T_i^+$ is $1+C$ for every vertex $i$. However, it is known that  $E_i T_i^+ = \frac{2m}{deg(i)}$, and therefore $deg(i)$ is constant, i.e. $G$ is regular. 

$(ii)$ Suppose that $G$ is highly symmetric and has a resistance-regular cut-vertex $i$. Let $a=E_iT_j = E_jT_i$ for any $j \sim i$ (our assumptions imply that this value is independent of the choice of $j$). Place an equivalence relation on the neighbors of $i$, where two vertices $j, k$ are equivalent if there is a path from $j$ to $k$ which does not pass through $i$. Since $i$ is a cut-vertex there are at least two equivalence classes, and at least one of them contains no more than $d/2$ vertices, where $d=deg(i)$. Choose $j$ in such an equivalence class. If we start a random walk at $i$ and let it go 1 step, it has at least a $1/2$ chance of moving to a vertex $k$ in an equivalence class not containing $j$, and in that case the only path from $k$ to $j$ must pass through $i$. Thus, $E_kT_j = 2a$, and it follows that $E_iT_j = a \geq \frac{1}{2}(1+2a) = a+1/2$, a contradiction.

$(iii)$ Suppose that $G$ is a highly symmetric graph with a cut-edge $(i,j)$. If we remove $(i,j)$ the graph now has two components, and we will call these components $G_i$ and $G_j$. $G$ has $m$ edges, and let $m_i, m_j$ be the number of edges in $G_i$ and $G_j$, respectively. We know that $R_{ij} = 1$, and since $G$ is highly symmetric we have $E_iT_j = m$. Now, if our random walk starts at $i$ and takes its first step to $j$ then $T_j = 1$, however if this does not happen then the first step is chosen randomly in $G_i$, and the random walk then lives in $G_i$ until such time as it returns to $i$. From that point on, it will take again an average of $m$ steps to hit $j$. If we set $d_i = deg(i)$, then we arrive at the equation

\begin{equation} \label{sleep}
    m = \frac{1}{d_i}(1) + \frac{d_i-1}{d_i}(E_iT^+_i + m), 
\end{equation}

where $T^+_i$ denotes the return time of a random walk in the graph $G_i$. Now, we know that $E_iT^+_i = \frac{2m_i}{d_i-1}$, and putting this into (\ref{sleep}) and rearranging yields $m = 1 + 2m_i$. Naturally the same holds for $m_j$, and this proves $(iii)$.

$(iv)$ Suppose now that $G$ has another cut-edge $(u,v)$. This edge must lie in either $G_i$ or $G_j$, so let us assume it lies in $G_i$. Then one of the components that remains when we remove $(u,v)$ must contain $G_j$ as well as the edge $(i,j)$, and therefore must contain more edges than the other component, contradicting $(iii)$. $\bullet$

There are a number of open questions regarding highly symmetric graphs, and in our opinion they are of considerable interest. Several are raised in \cite{george1}, and we repeat them here.

\begin{itemize}
    \item Does there exist a highly symmetric graph which is not walk-regular?
    
    \item Does there exist a highly symmetric graph which is not regular?
\end{itemize}

Georgakopoulos guessed that the answer to each of these questions is yes, but our guess is yes and no, respectively. A highly symmetric graph which is not regular would be a fascinating object indeed. A possible place to look for examples for the first question is in the family of edge-transitive regular graphs which are not vertex-transitive, since these are known to exist and be highly symmetric, but to our knowledge are not known to necessarily be walk-regular. Furthermore, in light of $(ii)-(iv)$ above, the following seem like natural questions to us.

\begin{itemize}
    \item Does there exist a highly symmetric graph with a cut-vertex?
    
    \item Does there exist a highly symmetric graph with a cut-edge (other than the trivial case of a tree on two vertices)?
\end{itemize}

Our guess on both of these is no, but we do not know how to prove it.

\end{document}